\documentclass{amsproc}

\usepackage{amsmath, amssymb, amsfonts}

\newtheorem{theorem}{Theorem}
\newtheorem{prop}[theorem]{Proposition}
\newtheorem{lemma}[theorem]{Lemma}
\newtheorem{remark}[theorem]{Remark}
\newtheorem{cor}[theorem]{Corollary}
\theoremstyle{definition}
\newtheorem{definition}[theorem]{Definition}
\theoremstyle{remark}

\def\bP{\mathbb P}
\def\bC{{\mathbb C}}

\def\C{{\mathcal C}}
\def\M{{\mathcal M}}
\def\H{\mathcal H}
\def\X{\mathcal X}
\def\sem{{\rtimes}}
\def\iso{{\, \cong\, }}

\def\<{\langle}
\def\>{\rangle}
\def\s{\sigma}
\def\l{\lambda}

\def\g{\gamma}
\def\p{\mathfrak p}
\def\sem{{\rtimes}}
\def\p{\mathfrak p}

\def\xs{{\Bbb o}}

\def\N{\mathcal N}

\def\a{\alpha}

\def\D{\Delta}

\begin{document}

\title{On the generic curve of genus 3}

\author{Tanush Shaska and Jennifer L. Thompson}

\address{300 Brink Hall, Department of Mathematics, University of Idaho, Moscow, ID, 83844.}
\email{tshaska@uidaho.edu}

\address{300 Brink Hall, Department of Mathematics, University of Idaho, Moscow, ID, 83844.}
%\email{jthompson@uidaho.edu}

%    General info
%\date{January 1, 1994 and, in revised form, June 22, 1994.}

\dedicatory{To Rachel and Adrianna.}


\begin{abstract} We study genus $g$ coverings of full moduli dimension of degree $d=[\frac {g+3} 2]$. There is a
homomorphism between the corresponding Hurwitz space $\H$ of such covers to the moduli space $\M_g$ of genus $g$
curves. In the case $g=3$, using the signature of such covering we provide an equation for the generic ternary
quartic. Further, we discuss the degenerate subloci of the corresponding Hurwitz space of such covers from the
computational group theory viewpoint. In the last section, we show that one of these degenerate loci corresponds
to the locus of curves with automorphism group $C_3$. We give necessary conditions in terms of covariants of
ternary quartics for a genus 3 curve to belong to this locus.
\end{abstract}

\maketitle

%\tableofcontents

\section{Introduction}

In this brief note we study non-hyperelliptic curves of genus 3. The main idea is to determine an equation for the
generic genus 3 curve starting from a covering of full moduli dimension.

Determining the monodromy group of a generic genus $g$ curve covering $\mathbb P^1$ is a problem with a long
history which goes back to Zariski and relates to Brill-Nother theory. Let $\X_g$ be generic curve of genus
$g$ and $f: \X_g \to {\mathbb P}^1$ a degree $n$ cover. Denote by $G:=Mon (f)$, the monodromy group of $f:
\X_g \to {\mathbb P}^1$. Zariski showed that for $g > 6$, $G$ is not solvable. For $ g \leq 6$ the situation
is more technical. This has been studied by many authors e.g., Fried, Guralnick, Neubauer, Magaard,
V\"olklein et al. The problem is open for $g=2$. The main question is to determine all possible  signatures
of the cover $f: \X_g \to {\mathbb P}^1$, where $\X_g$ is a generic curve. Covers with such signature are
called covers of full moduli dimension. In section 2, we provide a brief description of the main situation
and give at least one signature for each $g > 2$ which gives a cover of full moduli dimension. The degree of
the corresponding cover is the well known gonality of the curve. We show that the locus of such curves in
$\M_g$ does not intersect the hyperelliptic locus (cf. Proposition 4).

In section 3, we use this cover of full moduli dimension for $g=3$ to find an equation of the generic curve
given by
$$Y^3(X+a)+ Y^2(b X+c)+ Y(d X^2 +e X) +X^3 +f X^2 +X=0.$$ It is not clear whether the tuple $(a,b,c,d,e,f)$
is determined uniquely once we fix the moduli point. In other words, we don't know the degree of the extension
$k(a,b,c,d,e,f)/k (\p)$ where $k (\p)$ is the field of moduli. Investigating this question would give some insights
on the problem of field of moduli versus field of definition which is classical in algebraic geometry. We prove
computationally that $[k(a,b,c,d,e,f) : k (\p)] \le 3$. Hence the obstruction has at most degree 3.  We define six
$GL_3 (\bC)$-invariants $i_1, \dots , i_6$. However, it is not clear whether these invariants generate the field of
invariants of $GL_3 (\bC)$ on the space of ternary quartics. In other words, for a given curve it is not known if
$[k (\p): k(i_1, i_2, i_3, i_4, i_5, i_6)] =1$. If one were to show computationally that $[k(a,b,c,d,e,f) : k(i_1,
i_2, i_3, i_4, i_5, i_6)] =1$ then this would provide a proof that $(i_1, i_2, i_3, i_4, i_5, i_6)$ determines the
moduli point in $\M_3$.

In section 4, we discuss briefly the Hurwitz space of these coverings of full moduli dimension. We show that
this is an irreducible locus in $\M_3$ and compute the number of Nielsen classes. We show that there are 4
degenerate ways  the branch points of the covering can collapse. In each case we compute the number of
Nielsen classes. In the first three cases the monodromy group is $S_3$ and in the last case is the cyclic
group $C_3$.  The covering obtained in this case is a Galois covering and the curve has automorphism group
$C_3$.

The generic curve of genus 3 has no automorphisms.  In the last part we give a brief description of genus 3
non-hyperelliptic curves with automorphisms.  We make some remarks on the invariants of these curves and the
field of moduli and field of definition.

\bigskip

\textbf{Notation:}  By a curve we mean  an irreducible, smooth, projective curve  defined over an algebraically
closed field $k$. $\M_g$ denotes the  moduli space of smooth curves of genus $g$. We denote by $[\X_g]$ or by $\p$
the isomorphism class of $\X_g$, i.e. the corresponding point in $\M_g$. The cyclic group of order $n$ is denoted by
$C_n$ and the Klein 4-group by $V_4$.

\section{Covers of full moduli dimension}

Let $\M_g$ be the moduli space of curves of genus $g \geq 2$ and $\bP^1=\bP^1(\bC)$ the Riemann sphere. Let $\X$ be
a curve of genus $g$ and $\phi: \X \to \bP^1$ be a degree $n$ covering with $r$ branch points. By covering space
theory, there is a tuple $(\s_1, \dots , \s_r)$ in $S_n$ such that $\s_1 \cdot \cdot \cdot \s_r =1$ and $G:=\< \s_1,
\dots , \s_r\> $ is a transitive group in $S_n$. We call such a tuple the signature of $\phi$. Conversely, let $\s:
= (\s_1, \dots , \s_r)$ be a tuple in $S_n$ such that $\s_1 \cdot \cdot \cdot \s_r =1$ and $G:=\< \s_1, \dots ,
\s_r\> $ is a transitive group in $S_n$. We say a cover $\phi: \X \to \bP^1$ of degree $n$ is of type $\s$ if it has
$\s$ as signature. The genus $g$ of $\X$ depends only on $\s$ (Riemann-Hurwitz formula). Let $\H_\s$ be the set of
pairs $([f], (p_1, \dots , p_r)$, where $[f]$ is an equivalence class of covers of type $\s$, and $p_1, \dots , p_r$
is an ordering of the branch points of $\phi$.

By \cite{Be}, the \textbf{Hurwitz space}  $\H_\s$   carries a structure of quasi-projective variety (over
$\bC$). We have a morphism $$\Phi_\s: \H_\s \to \M_g$$ mapping $([f], (p_1, \dots , p_r)$ to the class $[\X]$
in the moduli space $\M_g$.

Each component of $\H_\s$ has the same image in $\M_g$. As in \cite{FMV}, we define  \textbf{moduli dimension
of $\s$} (denoted by $\dim (\s)$) as the dimension of $\Phi_\s(\H_\s)$; i.e., the dimension of the locus of
genus g curves admitting a cover to $\bP^1$ of type $\s$. We say $\s$ has   \textbf{full moduli dimension} if
$$\dim(\s)=\dim \M_g$$
We would like to explore the map $\Phi_\s$ for small $g$. The following two problems are significant.

\medskip

\noindent   \textit{Problem 1:}  Given a signature $\s$. Compute the moduli dimension of $\s$.

\smallskip

\noindent   \textit{Problem 2:} Given $g \geq 3$ and some $\s$ which has full moduli dimension. What can you
say about $\Phi_\s$?

\medskip

Fix $g\geq 3$. It is well known that the dimension of the moduli space $\M_g$ is $3g-3$. We want to find $\s$ of
full moduli dimension. It is known that $dim (\s) = r-3 $, when the quotient space has genus 0; see \cite{FMV}.
Thus, if a cover has full moduli dimension then it has $r=3g$ branch points. By $[x]$ we denote the integer part
of $x$. Then we have the following:

\begin{lemma} For any $g\geq 3$ there is a degree $d=[\frac {g+3} 2]$ cover
$$\psi_g: \X_g \to \bP^1$$
of full moduli dimension from a genus $g$ curve $\X_g$ such that it has $r=3g$ branch points and signature:

i) If $g$ is odd, then $\, \, \s=(\s_1, \dots , \s_r)$ such that $\s_1, \dots , \s_{r-1}\in S_d$ are
transpositions and $\s_r \in S_d$ is a 3-cycle.

ii) If $g$ is even, then $\, \, \s=(\s_1, \dots , \s_r)$ such that $\s_1, \dots , \s_r \in S_d $ are
transpositions.
\end{lemma}

\begin{remark} Note that the degree $d=[\frac {g+3} 2]$ is the minimum degree of a map
from a generic curve of genus $g$ to $\bP^1$, see \cite{Mu}. This is normally known as the {\it gonality}
$gon(\X_g)$ of a generic curve $\X_g$ of genus $g$.
\end{remark}

\begin{definition}
For a fixed $g \geq 3$ we call the cover $\psi_g: \X_g \to \bP^1$ as above the {\bf W-cover} associated to $g$.
\end{definition}

\begin{prop}
Let $\X_g$ be a hyperelliptic curve of odd genus $g\geq 3$. Then, there is no W-covering $\psi: \X_g \to \bP^1$.
\end{prop}

\proof Let $w$ be the hyperelliptic involution of $\X_g$. We assume that there exists a W-covering $\psi_g: \X_g
\to \bP^1$. Then there is $P \in \X_g$ such that the ramification index $e_\psi (P)$ of $P$ under $\psi$ is 3.
Denote the corresponding function fields by $K$ and $k(z)$. If $w$ does not fix $k(z)$, then $s=w(z)$ and
$K=k(z,w)$. Thus, $s$ and $w$ satisfy a symmetric polynomial $g(s,t)=0$ which has degree 3 in both $s$ and degree
2 in $w$. Then the curve $g(s,w)$ has genus $[ \frac {d-1} 2 ] \neq g$. Assume that $w$ fixes $k(x)$. Let $v$
denote its restriction in $k(x)$. Then, there exists a degree $d$ covering $\phi: \bP^1 \to \bP^1$ such that the
following diagram is commutative.
$$\begin{matrix} \X_g & \buildrel{w}\over{\to} & \bP^1 \\ \psi \, \, \downarrow & & \downarrow \, \phi \\
k(z) & \buildrel{v}\over{\to} & \bP^1
\end{matrix}$$
We denote by $a_0=v \circ \psi (P)= \phi \circ w (P)$. Then, in $(v \circ \psi)^{-1} (a_0)$ there could be
one point of index 6, one point of index 3 and three unramified points, or one point of index 3 and one point
of index 2. In $(\phi \circ w)^{-1} (a_0)$ there could be one point of index 6 or two points of index 3.
Thus, $e_{\phi \circ w} (P)=6$. Hence, $e_w (P)=2$ and $P$ is a Weierstrass point. Since $\X_g$ is
hyperelliptic then the Weierstrass gap sequence is $$1= n_1 < 3 < \dots < 2g-1.$$
Thus, there are no functions with a pole of order three at a Weierstrass point. Hence, there is no W-cover
for $\X_g$. \qed

%*******************************************************
\section{The case of genus 3}

%\subsection{Genus 3 curves and their invariants}

For the rest of this paper we will focus on the case $g=3$. Let $\C$ be a genus 3 curve defined over $k=\bC$ and $K$
its function field. It is well known that if $\C$ is a hyperelliptic curve then it corresponds to a binary octavic,
otherwise to a ternary quartic. Thus, the isomorphic classes of genus 3 curves correspond to (projective)
equivalence classes of binary octavics or ternary quartics. Let $S(n,r)$ denote the graded ring of (projective)
invariants of homogeneous polynomials of order $r$ in $n$ variables with coefficients in $\bC$. In this section we
will describe briefly  $S(3,4) $ (ternary quartics).

We denote by  $V$ be the set of ternary quartics. $V$ is a complex vector space of dimension 15. Let $\bC [V]$ be
the algebra of complex polynomial functions on $V$ and by $\bC [V]_n$ the set of homogeneous elements of degree
$n$ in $\bC [V]$. Then, $\bC [V] = \oplus_{n \ge 0} \, \bC[V]_n$ is a graded algebra.

The group $SL_3 (\bC)$ acts in a natural way in $V$. The invariant space under this action is the algebra
$\mathcal A$ of projective invariants of quartic plane curves.  $\mathcal A$ admits a homogeneous system of
parameters of degree 3, 6, 9, 12, 15, 18, 27. We denote these system of parameters by  $I_3, I_6, I_9, I_{12},
I_{15}, I_{18}, I_{27}$, where $I_j$'s are defined  as in Dixmier \cite{Di}.  They are homogenous polynomials of
degree 3, 6, 9, 12, 15, 18, 27 respectively and generate the ring of invariants $SL (3, 4) $ of ternary quartics.
For a general ternary quartic given by
\begin{small}
\begin{equation*}
\begin{split}
g(x,y, z) & = a_1 x^4 + 4 a_2 x^3y + 6 a_3 x^2 y^2 + 4  a_4 xy^3
+ a_5 y^4 + 4 a_6 x^3z+ 12a_7 x^2 y z + 12a_8 x y^2 z \\
& + 4a_9 y^3 z + 6a_{10} x^2 z^2 + 12a_{11} x y z^2
+6 a_{12} y^2 z^2 + 4 a_{13} x z^3 + 4 a_{14} y z^3 + a_{15}
z^4.
\end{split}
\end{equation*}
\end{small}
$I_3$ and $I_6$ are defined as follows:
\begin{equation*}
\begin{split}
I_3 (g) & = a_1 a_5 a_{15} + 3 (a_1 a_{12}^2 + a_5 a_{10}^2 a_{15} a_3^2) + 4 (a_2 a_9  a_{13} + a_6 a_4 a_{14}\\
& - 4 (a_1 a_9 a_{14} + a_5 a_6 \\
I_6 (g) & = \det (A), \textit{   where }
A  = \begin{bmatrix}
a_1    & a_3    & a_{10} & a_7    & a_6    & a_2    \\
a_3    & a_5    & a_{12} & a_9    & a_8    & a_4    \\
a_{10} & a_{12} & a_{15} & a_{14} & a_{13} & a_{11} \\
a_7    & a_9    & a_{14} & a_{12} & a_{11} & a_8    \\
a_6    & a_8    & a_{13} & a_{11} & a_{10} & a_7    \\
a_2    & a_4    & a_{11} & a_8    & a_7    & a_3    \\
\end{bmatrix} \\
\end{split}
\end{equation*}
see Dixmier \cite{Di} for definitions of $I_9, \dots , I_{27}$.

Two algebraic curves $\C$ and $\C'$ are isomorphic if there exists an $\alpha \in GL_3 ( \bC)$ such that $\C \iso
\C^\alpha$. We define $$i_1:= \frac {I_6} {I_3^2}, \quad i_2 := \frac {I_9} {I_3^3},  \quad i_3 := \frac {I_{12}}
{I_3^4}, \quad i_4 := \frac {I_{15}} {I_3^5}, \quad i_5 := \frac {I_{18}} {I_3^6}, \quad i_6 := \frac {I_{27}}
{I_3^9}.$$ Since $I_j$'s are $SL_3( \bC)$-invariants then $i_1, \dots , i_6$ are invariants under the $GL_3
(\bC)$-action. If $C \iso \C'$, then $i_j (\C) = i_j (\C')$, for $j= 1, \dots , 6$.

\begin{remark}
It is a difficult task of computing explicitly the invariants $I_j$ of a ternary quartic. This is because
their expressions are rather large. We will see that using the canonical form found in the following lemma,
this is much easier.
\end{remark}

%********************************************************************
\subsection{The generic curve of genus 3}

In this section we use the structure of W-covers to obtain an equation for the generic curve of genus 3. This
equation is given with 6 parameters (the moduli dimension is also 6), hence the name generic.

For $g=3$ a W-cover has degree three and 9 branch points. The signature is $\s=(\s_1, \dots , \s_9)$ where $\s_i\in
S_3$ is an transposition for $i=1, \dots , 8$ and $\s_9$ is the 3-cycle.

%Throughout this section $\C$ denotes a non-hyperelliptic genus 3 curve defined over $k=\bC$ and $K$ its function
%field.

\def\char{\mbox char \, }

\begin{lemma}
Let $\C$ be a generic curve of genus 3 defined over a field $L$ such that $\char L \neq 2 , 3$. Then, there is a
degree 3 covering $\psi: \C \to \bP^1$ of full moduli dimension. Moreover, $\C$ is isomorphic to a curve with affine
equation
$$Y^3(X+a)+ Y^2(b X+c)+ Y(d X^2 +e X) +X^3 +f X^2 +X=0$$
for $a,b,c,d,e,f \in L$ such that $\Delta \neq 0$, where $\Delta$ is the discriminant of the quartic.
\end{lemma}

\proof Let $P$ be a Weierstrass point on $\C$ and $K$ the function field of $\C$. Then exists a meromorphic function
$x$ which has $P$ as a triple pole and no other poles. Thus, $[K : L(x)]=3$. Consider $x$ as a mapping of $\C$ to
the Riemann sphere. We call this mapping $\psi: \C \to \bP^1$ and let $\infty$ be $\psi(P)$. From the
Riemann-Hurwitz formula we have that $\psi$ has at most 8 other branch points. There is also a meromorphic function
$y$ which has $P$ as a triple pole and no other poles. Thus the equation of $K$ is given by
\begin{equation} \label{eq_g3} F(x,y):=\g_1(x)\, y^3 + \g_2(x)\, y^2 + \g_3(x)\, y +\g_0 (x)=0
\end{equation}
where $\g_0 (x), \dots , \g_3 (x) \in L[x]$ and $deg (\g_i)=i$ for $i=1,2,3$, $deg(\g_0)=3$. The discriminant of
$F(x,y)$ with respect to $y$
$$D(F,y):=-27\, (\g_1 \, \g_0)^2 + 18\, \g_0\, \g_1\, \g_2\, \g_3 + (\g_2\, \g_3)^2 - 4\, \g_0 \, \g_2^3
- 4\, \g_1\, \g_3^3, $$
must have at most degree 8 since its roots are the branch points of $\psi:\C \to \bP^1$. Thus, we have
$$deg \, (\g_3 \g_2) \leq 4, \quad deg\, (\g_0 \g_2^3) \leq 8, \quad deg \, (\g_3^3\, \g_1)\leq 8. $$
If $deg\, (\g_2)=2$ then $deg \, (\g_3)\leq 2$ and $deg \, (\g_0)=0$. Thus, $deg\, (F,x)=2$. Then, $F(x,y)=0$
is not the equation of an genus 3 curve. Hence, $deg \, (\g_2) \leq 1$. Clearly, $deg\, (\g_3)\leq 2$. We
denote:
\begin{equation}
\begin{split} \g_1 (x):= & \, ax+b, \quad \quad \quad \quad \g_2 (x):= \, cx+d\\
\g_3 (x):= & \, ex^2+fx+g, \quad \g_0 (x):= \, hx^3+kx^2+lx+m\\
\end{split}
\end{equation}
Then we have $$F(x,y)=(ax+b)y^3+(cx+d)y^2+(ex^2+fx+g)y+ (hx^3+kx^2+lx+m)=0$$ We can make $g=0$ and $m=0$ by
the transformation $$ x \to x + r, \quad y \to y + s$$ \noindent such that
\begin{equation}
\begin{split}
3 b s^2+ 3a s^2r + 2c sr + er^2+ 2d s +f r +g=0 \\
a s^3 r + b s^3 + h r^3 + c s^2 r + d s^2 + e r^2 s + k r^2 f sr + lr + gs +m=0\\
\end{split}
\end{equation}
Thus we have
$$F(x,y)=(ax+b)y^3+(cx+d)y^2+(ex^2+fx)y+ x^3+kx^2+x = 0.$$
We can make $a=1$ and $h=1$ by the transformation
$$x \to r\, x,  \quad  y \to s \, y, \textit{   such that  } \quad
r^3= \frac 1 h, \quad s^3= \frac 1 {ar}.$$ Then,
$$F(x,y)=(x+b)y^3+(cx+d)y^2+(ex^2+fx)y+x^3+kx^2+lx = 0.$$
This completes the proof of the lemma.

\qed

\begin{lemma}
Let $\C$ and $\C'$ be two non-hyperelliptic genus 3 projective curves defined over $\bC$. If $\C$ and $\C'$ are
isomorphic then exists $\l \in k^*$ such that $$I_j'=- \l^{\frac j 3}  \cdot I_j,$$ for $ j=3, 6, 9, ,12, 15, 18, 27
$.
\end{lemma}

\proof

Let $\C$ be a genus 3 curve given by $F(X, Y, Z)=0$,  where
$$F=Y^3(X+aZ)+ Y^2(b XZ+cZ^2)+ Y(d X^2 Z  +e X Z^2) +X^3Z
+f X^2 Z^2 +X Z^3$$ defined over a $\bC$. Using some computer algebra system we can compute $I_3, \dots , I_{27}$.
The expressions are very large so we display only $I_3$ and $I_6$.

\begin{equation*}
\begin{split}
I_3 & = \frac 1 {144} \cdot (2fb^2+2cd^2-dbe-6bc-6ef-6adf+9d+9ea ) \\
I_6  & = \frac 1 {2^{12}\cdot 3^6} \cdot  ( 24cfb^3-4c^2d^4-4b^4f^2+36cd^3-216c^2e -36c^2b^2-36f^2e^2\\
&+144f^3c-81a^2e^2-648cf+108dfe-81d^2+24af^2b^2d+24cd^3fa-8cd^2fb^2\\
& +4ced^3b+4b^3dfe-144af^3b+24c^2d^2b+24b^2f^2e-b^2e^2d^2-36a^2f^2d^2\\
&+144c^2fd-108cbd-162aed-108ae^2f+72ce^2d-36b^2df+18bed^2+648abf\\
& -48cd^2fe+108ceab+18ae^2bd+144af^2ed-12cedb^2-36cead^2-12be^2df\\
& -36afeb^2+108a^2fed -72cfabd-12aed^2fb-108ad^2f  ) \\ \\
\end{split}
\end{equation*}

\noindent Let $\alpha  \in GL_3 (\bC)$, such that   $D:= det (\a)\neq 0$ and $\C'\iso \C^\a$. Then, $\C'$ has
equation $F (X', Y', Z')=0$ where
$$ \begin{pmatrix} X' \\ Y' \\ Z' \end{pmatrix}
= \a \, \begin{pmatrix} X \\ Y \\ Z \end{pmatrix}. $$
The  covariants $I_3, \dots , I_{27}$ of $F (X', Y', Z')$ are exactly
$$I_j'= D^4   \cdot I_j,$$ for $ j=3, 6, 9, ,12, 15, 18, 27 $.

\endproof

\begin{cor}
Let $\C$ and $\C'$ be two non-hyperelliptic genus 3 curves defined over $\bC$. If $\C$ and $\C'$ are isomorphic then
$$i_j (\C) = i_j (\C^\prime), \quad \textit{   for   } j=1, \dots , 6$$
\end{cor}

Since $i_1, \dots , i_6$ are $GL_2(\bC)$ invariants then $\bC(i_1, \dots , i_6) < \bC(a, b,c,d,e, f)$. The degrees
of $i_1, \dots , i_6$ are 6, 9, 12, 15, 18, and 27 respectively.  Since the degree
$[ \bC(i_1,i_2,i_3,i_4,i_5,i_6) : \bC(a,b,c,d,e,f)]$
must be a common factor of the above degrees, we have:

\begin{cor}
$ [ \bC(i_1,i_2,i_3,i_4,i_5,i_6) : \bC(a,b,c,d,e,f) ] \leq 3$
\end{cor}

One can attempt to find  $[\bC(a,b,c,d,e,f) : \bC(i_1,i_2,i_3,i_4,i_5,i_6)]$ computationally. However, we could not
accomplish this even using sophisticated computer algebra packages and a lot of computer power.  The upshot would be
to show that $\bC(a,b,c,d,e,f) = \bC(i_1,i_2,i_3,i_4,i_5,i_6)$, which we believe is true. This would have some very
important consequences. First, it would show that $i_1,i_2,i_3,i_4,i_5,i_6$ generate the field of invariants and
therefore describe a moduli point, which is still an open problem.  Second, for every moduli point $\p= (
i_1,i_2,i_3,i_4,i_5,i_6)$ it would provide a rational model of the curve over its field of moduli.

\begin{remark}
In this section we focused on curves defined over $\bC$ instead of a general field $k$. The main reason was that the
invariants of ternary quartics in \cite{Di} were defined over $\bC$.  One could make the necessary adjustments and
define Dixmier invariants over any field of characteristic $p \neq 2, 3$.
\end{remark}

\section{Hurwitz space of W-covers and moduli space $\M_3$}
%***************************************************************************
In this section we discuss the corresponding Hurwitz space and degenerations of the covering $\phi: \X_g \to
\bP^1$.

Let $\M_g$ be the moduli space of curves of genus $g \geq 2$ and $\bP^1=\bP^1(k)$ the Riemann sphere.  Let  $\phi:
\X_g \to \bP^1$ be a degree  $n$ covering with $r$ branch points. By covering space theory, there is a tuple
$(\s_1, \dots , \s_r)$ in $S_n$ such that $\s_1 \cdot \cdot \cdot \s_r =1$ and $G:=< \s_1, \dots , \s_r> $ is a
transitive group in $S_n$. We call such a tuple the {\it signature} of $\phi$. We say that a  permutation is of
type  $n^p$ if  it is a product  of $p$ disjoint  $n$-cycles.

Conversely, let $\s: = (\s_1, \dots , \s_r)$ be a tuple in $S_n$ such that $\s_1 \cdot \cdot \cdot \s_r =1$
and $G:=< \s_1, \dots , \s_r> $ is a transitive group in $S_n$. We say that a cover $\phi: \X \to \bP^1$ of
degree $n$ is of type $\s$ if it has $\s$ as signature. The genus $g$ of $\X$ depends only on $\s$
(Riemann-Hurwitz formula). Let $\H_\s$ be the set of pairs $([f], (p_1, \dots , p_r)$, where $[f]$ is an
equivalence class of covers of type $\s$, and $p_1, \dots , p_r$ is an ordering of the branch points of
$\phi$. The \textbf{Hurwitz space} $\H_\s$ is a quasiprojective variety. We have a morphism
$$\Phi_\s: \H_\s \to \M_g$$
mapping $([f], (p_1, \dots , p_r)$ to the class $[\X]$ in the moduli space $\M_g$.  Each component of $\H_\s$ has
the same image in $\M_g$.

We denote  by $C:=(C_1, \dots , C_r)$, where $C_i$ is the conjugacy class of $\s_i$ in $G$.  The set of Nielsen
classes $\N(G, C)$ is
$$\N(G, C):=\{ (\s_1, \, \dots \, , \s_r) \, |  \, \,  \s_i \in C_i,
\, G=< \s_1, \dots , \s_r >, \, \,  \s_1 \cdots \s_r=1  \} .$$

Fix a base point $\l_0\in \bP^1\setminus S$ where $S$ is the set of branch points. Then $\pi_1(\bP^1\setminus S)$
is generated by homotopy classes of loops $\g_1, \dots , \g_r$. The braid group acts on $\N(G, C)$ as
$$[\g_i]: \quad (\s_1, \, \dots \, , \s_r) \to (\s_1, \, \dots , \,
\s_{i-1},    \s_i \s_{i+1} \s_i^{-1}      , \s_i, \s_{i+2}, \dots , \s_r).$$

The orbits of this action are called the {\it braid orbits} and correspond to the irreducible components of $\H
(G,C):=\H_\s$.

\begin{lemma}
Let $C = ( 2^8, 3)$ and $\H_\s: = \H (S_3, C)$ be the corresponding Hurwitz space. Then, $\H_\s$ is an irreducible locus in $\M_3$
\end{lemma}

\proof The proof  is elementary and is based on the fact that there is one braid orbit of the braid
action. The computationally minded reader can use the braid program to compute the braid orbits,
see \cite{Ma}.
\endproof

\noindent{\bf Degenerate cases:} If some of the branch points of $\psi$ coalesce we have the following signatures:
$$
(3, 3, 2, 2, 2, 2, 2, 2), \quad (3, 3, 3, 2, 2, 2, 2), \quad (3, 3, 3, 3, 2, 2),  \quad (3, 3, 3, 3, 3).
$$
The information of each Hurwitz space is compiled in Table~\ref{tb1}.
\begin{table} [ht!]
\label{tb1}
\begin{center}
\renewcommand{\arraystretch}{1.24}
\begin{tabular}{||c|c|c|c||}
\hline \hline
Signature & Mon. group  & $\#$ Nielsen cl. & mod. dim.    \\
\hline \hline
 (3, $2^8$) & $S_3$ & 729 & 6 \\
 ($3^2, 2^5$) & $S_3$ & 36 &  5  \\
 ($3^3, 2^4$ ) & $S_3$ & 6 & 4  \\
 ($3^4, 2^2$) & $S_3$ & 6 & 3  \\
 ($3^5$) & $C_3$ & 5 &  2 \\
\hline \hline
\end{tabular}
\end{center}
\vspace{.5cm} \caption{The number of Nielsen classes in each Hurwitz space}
\end{table}

\bigskip

Proofs of the data in the Table~\ref{tb1} are simple exercises in combinatorial group theory. Moreover, one
can use the braid package in GAP written by K. Magaard.

\begin{remark}
These cases are studied in \cite{KK}. The main goal of that  paper is to study the Weierstrass points in each
case and provide equations for the curve. The equation of the curve in the last case is mistakenly identified.
\end{remark}

We find conditions on $a, b, c, d, e, f$ for each case. The discriminant of curve  given in Lemma 6 with
respect to $Y$ is a polynomial $\D (X)$ in $X$ given by
\begin{equation}
\begin{split}
\D (X) & = - X (27X^7 + A_6 X^6 + A_5 X^5 + A_4 X^4 + A_3 X^3 + A_2 X^2 + A_1 X + 4c^3)
\end{split}
\end{equation}
where $A_1, \dots A_6 $ are as follows:
\begin{small}
\begin{equation*}
\begin{split}
A_1 & =12bc^2+4c^3f-c^2e^2-18eca+27a^2     \\
A_2 & = 12b^2c+54a-18dca-2be^2c-18ec+54fa^2-18eba+4c^3-18ecfa\\
& +12bc^2f+4e^3a-2c^2de    \\
A_3 & = 27-18dcfa-18ebfa-18dc-18eb+108af+27f^2a^2-b^2e^2+12bc^2-\\
& c^2d^2-18eca-18dba -18ecf+12b^2fc+12de^2a+54a^2+4b^3+4e^3-4bdce    \\
A_4 & =-2bd^2c-18dcf+54f^2a+54fa^2-18ebf-18dbfa+12de^2+12d^2ae\\
& -18db-2b^2de+4b^3f +54f-18dca-18ec-18eba+108a+12b^2c     \\
A_5 & = 4d^3a+27a^2+54-18eb+108af-18dc-18dbf-b^2d^2+12d^2e\\
& -18dba+27f^2+4b^3    \\
A_6 & = 4d^3+54a-18db+54f    \\
\end{split}
\end{equation*}
\end{small}

The branch points of the cover $\phi : \C \to \bP^1$ coalesce when $\D (X)$ has multiple roots. Thus, its
discriminant $\D$ in X is $\D= 0$. There are four factors of the discriminant
$$\D = \D_1\cdot  \D_2\cdot  \D_3\cdot \D_4=0,$$
each corresponding to one of the degenerate cases in the previous table. We don't display them since they are easily
computed using Maple or any other computer algebra packages.

\subsection{Curves with automorphisms}
%&&&&&&&&&&&&&&&&&&&&&&&&&&&&&&&&&&&&&&&&&&&&&&&&&&&&&&&&&&&&&&&&&&&&

The generic curve has no automorphisms. For sake of completeness we will briefly describe non-hyperelliptic
genus 3 curves with automorphisms. Genus 3 hyperelliptic curves and their automorphisms are treated in
\cite{GS}.

There are several papers written on automorphism groups of genus 3 curves. The following table is taken from
\cite{KSSV} and classifies all such families.

\begin{table}[hbt!]
\begin{center}
\begin{tabular}{||c|c|c|c||}
\hline \hline
&&& \\
$G$  &  sig. & equation & Group \\
&&& ID \\

\hline  \hline

&&&\\

$V_4$ &    $(2^6)$ &$x^4+y^4+ax^2y^2+bx^2+cy^2+1=0$&(4,2) \\

$D_8$ &     $(2^5)$ & take\ $b=c$& (8,3)\\

$S_4$ &    $(2^3 ,3)$ & take\ $a=b=c$ &(24,12) \\

$C_4^2 \xs S_3$ &     $(2,3,8)$ &  \ take \, $a=b=c=0$ \, or\, $y^4=x(x^2-1)$ & (96,64) \\

\hline

$16$ &    $(2^3 ,4)$ &$y^4=x(x-1)(x-t)$&(16,13)        \\

$48$ &     $(2,3,12)$ &$y^4=x^3-1$ &(48,33)              \\

\hline

$C_3$    & $(3^5)$ &$y^3=x(x-1)(x-s)(x-t)$&(3,1)  \\

$C_6$     &    $(2,3,3,6)$ & take\ $s=1-t$& (6,2)     \\

$C_9$     & $(3,9,9)$ &$y^3=x(x^3-1)$&(9,1)\\

\hline &&&\\

$L_3(2)$ & $(2,3,7)$  &$x^3y+y^3z+z^3x=0$ &(168,42) \\

\hline

&&&\\
$S_3$ & $(2^4 ,3)$  & $a(x^4+y^4+z^4)+b(x^2y^2+x^2z^2+y^2z^2)+$&(6,1)\\

&&$c(x^2yz+y^2xz+z^2xy)=0$&\\
\hline &&&\\

$C_2$ & $(2^4)$ &$ x^4+x^2(y^2+az^2) + by^4+cy^3z+dy^2z^2$&(2,1)\\
&& $ +eyz^3+gz^4=0$, \  \ either $e=1$ or $g=1$ &\\

\hline \hline
\end{tabular}
\end{center}
%\vspace{.5cm}
\caption{Automorphism groups of genus 3 non-hyperelliptic curves}
\end{table}

In the table the cyclic group if order $n$ is denoted by $C_n$ and  $V_4$ denotes the Klein 4-group.  Each
group is identified also with the Gap identity number, see \cite{KSSV} for details. Each of the above cases
is an irreducible locus in $\M_3$ whose equation in terms of $i_1, \dots , i_6$ can be determined. Such
equations are for most cases large and we don't display them. The followings remarks are  rather easy to
check computationally.

\begin{remark}
If the automorphism group of a non-hyperelliptic genus 3 curve is isomorphic to $V_4$ and equation as in the above
table, then
\begin{equation}
\begin{split}
I_3 & = \frac 1 {36} \, (36+3c^2+3b^2+3a^2+abc)\\
I_6  & = \frac {abc}  {2^5 \, 3^6}   (108 - 3(a^2+b^2+c^2) +abc)
\end{split}
\end{equation}
The following three cases $D_8, S_4, C_4^2 \sem S_3$ are easily obtained.
\end{remark}

\begin{remark}
If the automorphism group of a non-hyperelliptic genus 3 curve is isomorphic to one of the following (16, 3), (48,
33), (3, 1), (6, 2), and (9, 1) then $I_3=I_6=0$. If the automorphism group is isomorphic to the group with Gap
identity (96, 64) then $I_3=1$ and $I_6=0$.
\end{remark}

\begin{remark}
The case when the automorphism group is isomorphic to $L_3 (2)$ corresponds to the Hurwitz curve (i.e., obtains
the maximum number of automorphisms).  In this case, $i_1 = - \frac 1 {16}$.
\end{remark}

\begin{remark}
In the last degenerate case, the covering $\phi: \X \to \bP^1$ with signature $$\s=(3, 3, 3, 3, 3)$$ is a Galois
covering. The curve $\X$ has automorphism group isomorphic to $C_3$. The Hurwitz space corresponding to these
coverings is the locus in $\M_3$ of genus 3 non-hyperelliptic curves with automorphism group $C_3$. These curves
have equation $$y^4= x(x-1)(x-s)(x-t). $$ The equation of this family is misidentified in \cite{KK}. The
covariants of such curves are $$I_3=I_6=0.$$
\end{remark}

It is an interesting question to investigate the field of moduli of these curves and see if that is a field of
definition. This is an open problem with a long history. In general the generic curve has no automorphisms and by
Weil's criterion its field of moduli is a field od definition.  Hence the curves on the above table are of special
interest. For the field of moduli of genus 3 hyperelliptic curves see \cite{GS}, \cite{Sh9}.

\bibliographystyle{amsalpha}

\end{document}